\documentstyle[11pt,twoside]{article}

\input{mssymb}
\makeatother

\newcommand{\mysec}[2]{%
\section*{\normalsize\hfil\sc {#1}. {#2}\hfill}%
\setcounter{theo}{0}\setcounter{equation}{0}\setcounter{section}{#1}%
\typeout{#1. #2}
\noindent}

\newcommand{\Proof}{\par\noindent{\em Proof. }}
\newcommand{\eop}{\nopagebreak\hspace*{\fill}$\Box$}

\newtheorem{theo}{Theorem}[section]
\newtheorem{lemma}[theo]{Lemma}
\newtheorem{cor}[theo]{Corollary}

\newtheorem{prop}[theo]{Proposition}
\newtheorem{definition}[theo]{Definition}
\newenvironment{defi}{   
   \begin{definition}
   \begin{em}}{
   \end{em}
   \end{definition}}

\newcounter{abc}   
\newcounter{iiiii} 

\newenvironment{aequivalenz}
{\setcounter{iiiii}{0}
\begin{list}%
{{\rm (\roman{iiiii})}}
{\usecounter{iiiii}
\parsep=0pt plus 1pt
\topsep=1pt plus 2pt minus 1pt
\itemsep=1pt plus 2pt minus 1pt
\leftmargin=3\baselineskip
\labelsep=.6\baselineskip
\labelwidth=2.4\baselineskip
\rightmargin 0pt}%
}
{\end{list}}

\newenvironment{statements}%
{\setcounter{abc}{0}
\begin{list}%
{{\rm (\alph{abc})}}
{\usecounter{abc}
\parsep=0pt plus 1pt
\topsep=1pt plus 2pt minus 1pt
\itemsep=1pt plus 2pt minus 1pt
\leftmargin=3\baselineskip
\labelsep=.6\baselineskip
\labelwidth=2.4\baselineskip
\rightmargin 0pt}%
}
{\end{list}}

{\begin{list}%
{$\bullet$}
{\leftmargin=3\baselineskip
\labelsep=\baselineskip
\labelwidth=2.5\baselineskip
\parsep=0pt plus 1pt
\topsep=1pt plus 2pt minus 1pt
\itemsep=1pt plus 2pt minus 1pt
\rightmargin 0pt}%
}
{\end{list}}

\newif\ifrefsc
\let\thebibliographyalt=\thebibliography                                %
\def\thebibliography#1                                                  %
 {\ifrefsc
 \section*{\normalsize\hfil\sc References\hfill}
 \small
 \list                    %
 {[\arabic{enumi}]}{\settowidth\labelwidth{[#1]}\leftmargin\labelwidth  %
 \advance\leftmargin\labelsep                                           %
 \usecounter{enumi}}                                                    %
 \def\newblock{\hskip .11em plus .33em minus .07em}                     %
 \sloppy\clubpenalty4000\widowpenalty4000\vfuzz 2pt                      %
 \sfcode`\.=1000\relax                                                  %
 \else\thebibliographyalt{#1}\fi}                                         %

\makeatletter
\def\ps@smallheadings{            
  \def\@oddfoot{}                 
  \def\@evenfoot{}                
  \def\@evenhead{\hbox to \textwidth {%
  \vbox{\hbox to \textwidth 
        {\thepage\hfil\strut{\footnotesize\sc\Kurzautor}\hfil}\vss}}}
  \def\@oddhead{\hbox to \textwidth {%
\vbox{\hbox to \textwidth 
        {\hfil{\footnotesize\sc\Kurztitel}\hfil\strut \thepage}\vss}}}
    }
\makeatother

\makeatletter

\def\eqalignno#1{\displ@y \tabskip\@centering
  \halign to\displaywidth{\hfil$\@lign\displaystyle{##}~$\tabskip\z@skip
    &$\@lign\displaystyle{{}##}$\hfil\tabskip\@centering
    &\llap{$\@lign##$}\tabskip\z@skip\crcr
    #1\crcr}}

\makeatother

\def\ersteSeite{\vspace*{32pt plus 2pt minus 2pt}\begin{center}
{\Large\sf\Titel}\\[15pt]{\sc\Autor}\\[26pt plus 2pt minus 2pt]
\end{center}\Abstrakt\vspace{0pt plus 2pt }\thispagestyle{empty}}

\def\Abstrakt{\begin{quote}\small\noindent{\sc Abstract.}
\Abstrakttext\end{quote}}

\newcommand{\N}{{\Bbb N}}
\newcommand{\R}{{\Bbb R}}

\newcommand{\Id}{{\rm Id}}

\newcommand{\eps}{\varepsilon}

\newcommand{\Sig}{\Sigma}

\newcommand{\loglike}[1]{\mathop{\rm #1}\nolimits}

\newcommand{\ex}{\loglike{ex}}

\newcommand{\co}{\loglike{co}}    
\newcommand{\coq}{\loglike{\overline{\co}}}
\newcommand{\lin}{\loglike{lin}}  
\newcommand{\linq}{\loglike{\overline{\rm lin}}}

\newcommand{\bea}{\begin{eqnarray*}}
\newcommand{\eea}{\end{eqnarray*}}
\newcommand{\beq}{\begin{equation}}
\newcommand{\eeq}{\end{equation}}
\newcommand{\begsta}{\begin{statements}}
\def\endsta{\end{statements}}
\newcommand{\begaeq}{\begin{aequivalenz}}
\def\endaeq{\end{aequivalenz}}

\newcommand{\dopu}{{:}\allowbreak\ }

\newcommand{\Norm}{|\mkern-2mu|\mkern-2mu|}

\newcommand{\rest}[2]{#1\raisebox{-0.3ex}{\mbox{$\mid_{#2}$}}} 

\def\DP{Daugavet property}
\newcommand{\Fin}{{\rm FIN}(\N)}

\refsctrue

\begin{document}


\def\Titel{ Banach spaces with the \DP }
\def\Kurztitel{ \Titel }
\def\Autor{ Vladimir~M.~Kadets\footnote{The work of this author was
done during his visit to Freie Universit\"at Berlin,
supported by a grant from the {\it Deutscher Akademischer
Austauschdienst}.}, 
Roman~V.~Shvidkoy, Gleb~G.~Sirotkin 
and~Dirk~Werner  }
\def\Kurzautor{ V.~M.~Kadets, R.~V.~Shvidkoy, G.~G.~Sirotkin and D.~Werner  }
\def\Abstrakttext{
A Banach space $X$ is said to have the Daugavet property if every
operator $T\dopu X\to X$ of rank~$1$ satisfies $\|\Id+T\| = 1+\|T\|$.
We show that then every weakly compact operator satisfies this
equation as well and that $X$ contains a copy of $\ell_{1}$. However,
$X$ need not contain a copy of $L_{1}$. We also study pairs of spaces
$X\subset Y$ and operators $T\dopu X\to Y$ satisfying
$\|J+T\|=1+\|T\|$, where $J\dopu X\to Y$ is the natural embedding.
This leads to the result that a Banach space with the Daugavet
property does not embed into a space with an unconditional basis. In
another direction, we investigate spaces where the set of operators
with $\|\Id+T\|=1+\|T\|$  is as small as possible and give
characterisations in terms of a smoothness condition.
}

\ersteSeite

\mysec{1}{Introduction}%
It is a remarkable result due to Daugavet \cite{Daug}  that the norm
identity
\beq\label{eq1.1}
\|\Id+T\|=1+\|T\|,
\eeq
which has become known as the Daugavet equation, holds for compact
operators on $C[0,1]$; shortly afterwards the same result for compact
operators on $L_{1}[0,1]$ was discovered by Lozanovskii \cite{Loz}.
Over the years, (\ref{eq1.1}) was extended to larger classes of operators
on various spaces (see for instance 
\cite{Abra2}, \cite{FoiSin}, \cite{Kadets}, \cite{PliPop},
\cite{Weis2}, \cite{Dirk10}, \cite{Woj92} and the references in these
papers for more information); in particular, the Daugavet equation
holds for operators not fixing a
copy of $C[0,1]$ defined on certain ``large'' subspaces of $C(K)$,
where $K$ is a compact space without isolated points,  and
for operators not fixing a copy of $L_{1}[0,1]$
defined on certain ``large'' subspaces of $L_{1}[0,1]$
(\cite{KadPop}, \cite{WeisDirk}).

In Section~2 of this paper we show that the validity of (\ref{eq1.1})
for weakly compact operators (as a matter of fact even for strong
Radon-Nikod\'ym operators) follows already from the corresponding
statement for operators of rank~$1$ which in turn can be verified by
means of a simple lemma. This has several interesting consequences and implies
for instance that $X$ contains a copy of $\ell_{1}$ if all the
rank-$1$ operators on $X$ satisfy the Daugavet equation. 
The investigation of the Daugavet equation for rank-$1$ operators is
suggested by results in Wojtaszczyk's paper \cite{Woj92}.
Actually, we
take a broader approach and consider a Banach space $Y$ together with
a closed subspace $X$ and the canonical embedding $J\dopu X\to Y$,
and we study the equation
\beq\label{eq1.2}
\|J+T\|=1+\|T\|
\eeq
for classes of operators from $X$ into $Y$. If (\ref{eq1.2}) holds
for operators of rank~$1$, we say that the pair $(X,Y)$ has the \DP.
We study heredity properties of Daugavet pairs and spaces and
prove that the \DP\ is inherited by $M$-ideals and by subspaces with
a separable annihilator.

To explain the relevance of studying (\ref{eq1.2}) we recall that, by
an argument from \cite{Kadets}, the well-known fact that neither
$C[0,1]$ nor $L_{1}[0,1]$ have unconditional bases can easily be
deduced from the \DP\ for finite-rank operators on these spaces. In
order to obtain the more general result that $C[0,1]$ and $L_{1}[0,1]$
do not even embed into spaces having an unconditional basis it is
necessary to investigate (\ref{eq1.2}). This was done in
\cite{KadShv}. Using techniques from that paper we now prove that if
$X$ is a subspace of a separable space $Y$ and $X$ has the \DP, then
$Y$ can be renormed so that the new norm coincides with the original
one on $X$ and the pair $(X,Y)$ has the \DP. This implies that no
space with the \DP\ embeds into a space with an unconditional basis.

Section~3 deals with quotient spaces of $L_{1}=L_{1}[0,1]$ by
``small'' subspaces. In particular, we exhibit an example of a space
with the \DP\ not containing a copy of $L_{1}$; our example is $L_{1}/Y$
with $Y$ a space constructed by Talagrand in his work on the three
space problem \cite{Tala5}.

In Section~4 we extend results from \cite{AbraAB} and \cite{Kadets} on
the following question. If $T\dopu X\to X$ is an operator with
$\|T\|\in\sigma(T)$, then $T$ is easily seen to satisfy
(\ref{eq1.1}). We say that $X$
has the anti-\DP\ (for a class of operators) if no other operators (in
this class) satisfy (\ref{eq1.1}). We show that the anti-\DP\ is
closely related to a smoothness property of $X$ introduced in
Definition~\ref{D3.1} below. In fact, we are able to characterize the
anti-\DP\ for compact operators along these lines.

We use standard notation such as 
$B_{X}$  and $S_{X}$ for the unit ball and the unit sphere of 
a Banach space $X$, and we employ the notation
$$
S(x^{*},\eps)= \{x\in B_{X}\dopu x^{*}(x)\ge 1-\eps\}
$$
for the slice of $B_{X}$ determined by $x^{*}\in S_{X^{*}}$ and $\eps>0$.
$\ex C$ stands for the set of extreme points of a set $C$.
In this paper we deal with real Banach spaces although our results
extend to the complex case with minor modifications.

\mysec{2}{The Daugavet property}%
Let $X$ be a subspace of a Banach space $Y$  and let $J\dopu X\to Y$
denote the inclusion operator. We say that the pair $(X,Y)$ has the
{\em Daugavet property\/} for a class $\cal M$ of operators, where ${\cal
M}\subset L(X,Y)$, if
\beq\label{eq1}
\|J+T\| = 1 + \|T\|
\eeq
for all $T\in{\cal M}$.
If $X=Y$, we simply say that $X$ has the Daugavet property with
respect to $\cal M$, and if $\cal M$ is the class of rank-$1$
operators, we just say that $X$ or $(X,Y)$ has the \DP.

It is well known \cite{AbraAB} that if $T$ satisfies (\ref{eq1}), then so
does $\lambda T$ for every $\lambda>0$. In particular, if $\cal M$ is
a cone, it is sufficient
to verify (\ref{eq1}) for all $T\in\cal M$ with $\|T\|=1$.
It is obvious that $X$ has the \DP\ once $X^{*}$ has it, but the
converse fails (take $X=C[0,1]$); also, a 1-complemented subspace of a
space with the \DP\ might fail it (consider one-dimensional
subspaces).

Here is the key lemma of our paper.

\begin{lemma}\label{L1}
Let $(X,Y)$ have the Daugavet property. Then:
\begsta
\item
For every $y_{0}\in S_{Y}$ and for every slice
$S(x^*_{0},\eps_{0})$ of $B_{X}$ there is another slice
$S(x^*_{1},\eps_{1}) \subset S(x^*_{0},\eps_{0})$ of $B_{X}$  such that for
every $x\in S(x^*_{1},\eps_{1})$  the inequality $\|x+y_{0}\| \ge
2-\eps_{0}$ holds.
\item
For every $x_{0}^{*}\in S_{X^{*}}$ and for every weak$^{\,*}$ slice
$S(y_{0},\eps_{0})$ of $B_{Y^{*}}$ $($where $y_{0}\in S_{Y}\subset
S_{Y^{**}})$ there is another weak$^{\,*}$ slice
$S(y_{1},\eps_{1}) \subset S(y_{0},\eps_{0})$ of $B_{Y^{*}}$  such that for
every $y^{*}\in S(y_{1},\eps_{1})$  the inequality $\|x_{0}^{*}+
\rest{y^{*}}{X}\| \ge 2-\eps_{0}$ holds.
\endsta
\end{lemma}

\Proof
Both parts are proved in a very similar fashion; so we only present
the proof of~(a).

Define $T\dopu X\to Y$ by $Tx=x^*_{0}(x)y_{0}$. Then $\|J^{*}+T^{*}\|
=\|J+T\|=2$, so there is a functional $y^*\in S_{Y^{*}}$ such that
$\|J^{*}y^* + T^{*}y^*\|\ge  2-\eps_{0}$ and $y^*(y_{0})\ge0$.
Put
$$
x^*_{1} = \frac{J^{*}y^* + T^{*}y^*}{\|J^{*}y^* + T^{*}y^*\|}, \quad
\eps_{1} = 1 - \frac{2-\eps_{0}}{\|J^{*}y^* + T^{*}y^*\|} .
$$
Then we have, given $x\in S(x^*_{1},\eps_{1})$,
$$
\langle (J^{*}+T^{*})y^*, x \rangle \ge  (1-\eps_{1}) \|J^{*}y^* + T^{*}y^*\|
= 2-\eps_{0},
$$
therefore
\beq\label{eq2}
y^*(x) + y^*(y_{0})x^*_{0}(x) \ge  2-\eps_{0},
\eeq
which implies that $x^*_{0}(x)\ge 1-\eps_{0}$, i.e., $x\in S(x^*_{0},
\eps_{0})$.
Moreover, by (\ref{eq2}) we have $y^*(x)+y^*(y_{0}) \ge  2-\eps_{0}$
and hence $\|x+y_{0}\| \ge  2-\eps_{0}$.
\eop

\bigskip
It is evident that the converse statement is valid, too. For future
reference we will record this in the following simplified version of
Lemma~\ref{L1}.

\begin{lemma}\label{L2}
The following assertions are equivalent:
\begaeq
\item
The pair $(X,Y)$ has the Daugavet property.
\item
For every $y\in S_{Y}$, $x^{*}\in S_{X^{*}}$ and $\eps>0$
there is some $x\in S_{X}$ such that $x^{*}(x)\ge 1-\eps$ and $\|x+y\|\ge
2-\eps$.
\item
For every $y\in S_{Y}$, $x^{*}\in S_{X^{*}}$ and $\eps>0$
there is some $y^*\in S_{Y^{*}}$ such that $y^*(y)\ge 1-\eps$ and
$\|x^{*}+\rest{y^*}{X}\|\ge 2-\eps$.
\endaeq
\end{lemma}

One consequence of these lemmas is that every slice of $B_{X}$ and
every weak$^{*}$ slice of $B_{X^{*}}$ has
diameter~$2$  if $(X,Y)$  has the \DP. In
particular, $X $ fails the Radon-Nikod\'ym property, a fact
originally due to Wojtaszczyk \cite{Woj92}. Likewise, $X^{*}$ fails
the Radon-Nikod\'ym property.
We shall return to this circle of ideas in Theorem~\ref{T5}.

We now come to our first main result.

\begin{theo}\label{T3}
If the pair $(X,Y)$ has the Daugavet property,
then $(X,Y)$ has the Daugavet property for weakly compact operators.
\end{theo}

Actually, the proof will show that the \DP\ holds for an even larger
class of operators, namely the strong Radon-Nikod\'ym operators,
meaning operators~$T$ for which $\overline{T(B_{X})}$ is a
Radon-Nikod\'ym set. 
It would be interesting to decide whether the result also extends to
operators not fixing a copy of $\ell_{1}$.

\bigskip
\Proof
Let $T\dopu X\to Y$ be a weakly compact operator with $\|T\|=1$. Then
$K=\overline{T(B_{X})}$ is weakly compact and therefore coincides
with the  closed
convex hull of its denting points; in fact, $K$ is the closed convex
hull of its strongly exposed points \cite{Bour-sex}. 
So for every $\eps>0$ there
is a denting point $y_{0}$ of $K$  with $\|y_{0}\|> \sup\{\|y\|\dopu
y\in K\} -\eps =1-\eps$, and for some $0<\delta<\eps$ 
there is a slice $S=\{y\in K\dopu
y^*(y) \ge  1-\delta \}$ of $K$ containing $y_{0}$ and having diameter
${<\eps}$; here $y^*\in Y^{*}$ and $\sup_{y\in K} y^*(y) =1$. 
Consider $x^{*}=T^{*}y^*$. By construction $\|x^{*}\|=1$ and 
\bea
T(S(x^{*},\delta)) &=&
\{Tx\dopu x\in B_{X},\ x^{*}(x)\ge 1-\delta \} \\
&=&
\{Tx \dopu x\in B_{X},\ y^*(Tx) \ge  1-\delta \} ~\subset~ S.
\eea
So for every $x\in S(x^{*},\delta)$ we have $\|Tx\|\ge 1-2\eps$.
Now by Lemma~\ref{L2} select an element $x_{0}\in S(x^{*},\delta)$ such
that $\bigl\| \,x_{0} +  y_{0}/\|y_{0}\|\, \bigr\| \ge  2-\delta$  and hence
$\|x_{0}+y_{0}\|\ge  2-2\eps$. But $Tx_{0}\in S$, so
$\|Tx_{0}-y_{0}\|<\eps$, and we have
$$
\|J+T\| \ge \|x_{0}+Tx_{0}\| \ge \|x_{0}+y_{0}\| - \eps \ge 2-3\eps.
\eqno\Box
$$

\smallskip\noindent
{\em Example.}
If $X$ is a Banach space and $(\Omega, \Sigma,\mu)$ is a non-atomic
measure space, then $Y:= L_{1}(\mu,X)$ has the \DP\ for weakly
compact operators.
This is a special case of a result due to Nazarenko \cite{Naz}; using our
preceding results we can prove this now in a few lines. Even in the
case of the scalar-valued function space $L_{1}(\mu)$, for which
other proofs have appeared for instance in \cite{BabPich},
\cite{Holub3} or \cite{Loz}, our argument is shorter.

In fact, let $y\in S_{Y}$ and $y^{*}\in S_{Y^{*}}$. The functional
$y^{*}$
can be represented by a weak$^{*}$ measurable function $\varphi$
taking values in $X^{*}$. For $\eps>0$, find a measurable subset $B$
of $\Omega$ such that $\|\chi_{B} y\|_{L_{1}} \le\eps/2$ and
$\|\chi_{B}\varphi\|_{L_{\infty}}\ge 1-\eps/2$, and pick $x\in S_{Y}$
so that $\chi_{B}x=x$ and $\langle \varphi,x \rangle \ge1-\eps$.
Since clearly $\|x+y\|\ge 2-\eps$, condition~(ii) of Lemma~\ref{L2} is
fulfilled.

By a similar argument, one can reprove the result from \cite{Kadets} that
$C(K,X)$ has the \DP\ if the compact space $K$ has no isolated
points.
\eop                                       

\bigskip
We now show that the \DP\ automatically extends to certain larger
range spaces.

If $K$ is a compact topological space, then we denote by
$\ell_{\infty}(K)$ the sup-normed space of bounded real-valued
functions on $K$ and by $m(K)$ the closed subspace consisting of
those $f\in \ell_{\infty}(K)$ for which $\{t\dopu f(t)\neq0\}$ is of
first category. Let us consider the quotient space
$$
m_{0}(K):= \ell_{\infty}(K)/m(K).
$$
It follows from the Baire category theorem that the restriction of
the quotient map $Q\dopu \ell_{\infty}(K)\to m_{0}(K)$ to $C(K)$ is
an isometry. It is proved in \cite{KadShv} that the pair $(C[0,1],
m_{0}(K))$, where $K$ is the unit ball of $C[0,1]^{*}$ in its
weak$^{*}$ topology, has the \DP, even for the class of so-called
narrow operators, and likewise for $L_{1}[0,1]$.

\begin{prop}\label{T6}
Let $(X,Y)$ have the \DP\ and denote by $K$  the weak$^{\,*}$ closure
of $\ex B_{Y^{*}}$ in $Y^{*}$. Then $(X, m_{0}(K))$ has the \DP, too.
\end{prop}

\Proof
The canonical map $J_{1}$ of $Y$ into $m_{0}(K)$ is the
composition $QJ_{0}$ where $J_{0}\dopu Y\to C(K)$ is defined by
$(J_{0}y)(y^{*})= y^{*}(y)$ and $Q$ is the (restriction of the) above
quotient map. As remarked above, $J_{1}$ is an isometry. To prove the
proposition we will employ Lemma~\ref{L2}. Thus, $x^{*}\in
S_{X^{*}}$, $\eps>0$ and $[f]\in S_{m_{0}(K)}$ are given; of course,
$[f]$ denotes the equivalence class of $f\in \ell_{\infty}(K)$ in
$m_{0}(K)$. We may assume, without loss of generality, that
$$
A=\{t\in K\dopu f(t)>1-\eps/2 \}
$$
is a set of second category. Now, there is an open set $V\subset K$
such that whenever $U\subset V$ is open, then $U\cap A$ is of second
category \cite[p.~202]{Kelley}.

Next we use the fact that an extreme point of a compact convex set has a
neighbourhood base consisting of slices; 
this follows from the converse to the Krein-Milman theorem 
(see \cite[p.~107]{Choq2} for details). 
Thus, $V$ contains a weak$^{*}$ slice; i.e., 
there are $y_{0}\in S_{Y}$ and $\eps_{1}<\eps/2$ such that 
\beq\label{eqm01}
w^{*}\in K,\ w^{*}(y_{0})>1-\eps_{1} \quad\Rightarrow\quad
w^{*}\in V .
\eeq
Since the pair $(X,Y)$ has the \DP, we get from Lemma~\ref{L2} some
$x_{0}\in S_{X}\cap S(x^{*},\eps/2)$ such that $\|x_{0}+y_{0}\|> 2-
\eps_{1}$; we wish to show that $\|[J_{0}x_{0}+f]\|\ge2-\eps$. 
In fact, for some $w_{0}^{*}\in K$ we
have  $w_{0}^{*}(x_{0})+ w_{0}^{*}(y_{0})
> 2-\eps_{1}$ and therefore $w_{0}^{*}(y_{0})>1-\eps_{1}$; thus
$w_{0}^{*}\in V$ by (\ref{eqm01}). Pick a neighbourhood $U\subset V$
of $w_{0}^{*}$ so that
$$
w^{*}(x_{0}) + w^{*}(y_{0}) > 2-\eps_{1} \qquad\forall w^{*}\in U.
$$
In particular $w^{*}(x_{0})>1-\eps_{1}$, and since $U\cap A$  is of
second category by construction of $V$, we have
\bea
\|[J_{0}x_{0}+f]\| &\ge&
\inf \{|w^{*}(x_{0})+f(w^{*})|\dopu w^{*}\in U\cap A\} \\
&\ge&
1-\eps_{1} + 1 - \frac\eps2 ~\ge~ 2-\eps,
\eea
and we are done.
\eop
\bigskip

Clearly, for every space between $Y$ and $m_{0}(K)$ the \DP\ holds,
as well; in particular, the pair $(X, C(\overline{\ex B_{Y^{*}}}))$
has the \DP. 
It remains open whether in general $(X,C(B_{Y^{*}}))$ has the \DP;
this is true for $X=Y=C[0,1]$ and $X=Y=L_{1}[0,1]$ \cite{KadShv}.

Proposition~\ref{T6} permits the following renorming theorem.

\begin{theo}\label{T6a}
Let $X\subset Y$ be separable Banach spaces, and suppose $X$ has the
\DP.
\begsta
\item
If $K= \overline{\ex B_{X^{*}}}$ and $m_{0}(K)$ is as above, then
there is an isomorphic embedding $S\dopu Y\to m_{0}(K)$ with
$\rest{S}{X}=\rest{Q}{X}$, the restriction of the quotient map from
$\ell_{\infty}(K)$ onto $m_{0}(K)$.
\item
The space $Y$ can be renormed so that the new norm coincides with the
original one on $X$ and $(X,Y)$ has the \DP\ for the new norm.
\endsta
\end{theo}

\Proof
Part~(a) is proved in \cite[Th.~3.1]{KadShv} for $K={B_{X^{*}}}$ (no
matter if $X$ has the \DP\ or not); but
the proof works for all symmetric isometrically norming weak$^{*}$
closed subsets of $B_{X^{*}}$
without isolated points. So, here it is enough to check that $\ex B_{X^{*}}$
does not contain isolated points.

Assume that $x_{0}^{*}$ is an isolated extreme point. Since $x_{0}^{*}$
has a neighbourhood base of slices, there are some $\xi\in S_{X}$ and
$\eps>0$ such that the only extreme point satisfying $x^{*}(\xi)\ge
1-\eps$ is $x^{*}=x_{0}^{*}$. Now the \DP\ yields an element $x\in
S_{X}$ with $x_{0}^{*}(x)\ge 1-\eps$  and $\|x-\xi\|\ge 2-\eps$. Pick
an extreme point $x^{*}$ so that $x^{*}(x-\xi)=\|x-\xi\|$; then we
have $(-x^{*})(\xi)\ge1-\eps$ and, consequently, $x^{*}=-x_{0}^{*}$
so that $1-\eps\le x^{*}(x)=-x_{0}^{*}(x)\le0$: a contradiction.

(b) By part~(a) and
 Proposition~\ref{T6}, $\Norm y \Norm = \|Sy\|$ defines a norm
with the required properties.
\eop
\bigskip

This theorem entails some information on the non-existence of
unconditional expansions. We first formulate a lemma that gives a
quantitative version of \cite[Lemma~3.6]{KadPop}.

\begin{lemma}\label{L6aa}
Let $X\subset Y$  be Banach spaces with $J\dopu X\to Y$ the natural
embedding. Suppose that the pair $(X,Y)$ has the \DP\ with respect to
a subspace ${\cal M}\subset L(X,Y)$  of
operators.  Let $T=\sum_{n=1}^{\infty} T_{n}$  be a pointwise
unconditionally convergent series of operators $T_{n}\in \cal M$.
Then $\|J+T\|\ge1$. In particular,
$J$ cannot be represented as a pointwise unconditionally convergent
sum $J=\sum_{n=1}^{\infty}T_{n}$ with $T_{n}\in \cal M$ for all
$n\in\N$.
\end{lemma}

\Proof
Denote by $\Fin$ the set of finite subsets of $\N$. By the 
Banach-Steinhaus theorem, the quantity
$$
\alpha = \sup \left\{ \biggl\| \sum_{n\in A} T_{n} \biggr\| \dopu
                    A\in \Fin \right\} 
$$
is finite, and whenever $B\subset \N$, then
$$
\biggl\| \sum_{n\in B} T_{n} \biggr\| 
\le
 \sup \left\{ \biggl\| \sum_{n\in A} T_{n} \biggr\| \dopu
                    A\in \Fin,\ A\subset B \right\} 
\le\alpha.
$$
Let $\eps>0$ and pick $A_{0}\in \Fin$ such that $\|\sum_{n\in A_{0}}
T_{n} \| \ge \alpha-\eps$. Then we obtain from the \DP
$$
\|J+T\| \ge
\biggl\| J+\sum_{n\in A_{0}} T_{n} \biggr\| -
\biggl\| \sum_{n\notin A_{0}} T_{n} \biggr\| 
\ge
1 + \biggl\| \sum_{n\in A_{0}} T_{n} \biggr\| -
\alpha 
\ge
1-\eps,
$$
which proves the lemma.
\eop
\bigskip

We could not decide whether in general $\|J+T\|=1+\|T\|$ in
Lemma~\ref{L6aa}. This is the case if the expansion is
1-unconditional, because for $\lambda_{n}=-1$ if $n\in A_{0}$ and
$\lambda_{n}=1$ otherwise
\bea
\|J+T\| &\ge& 1 + 2 \biggl\| \sum_{n\in A_{0}} T_{n} \biggr\| -
\biggl\| \sum_{n\in\N} \lambda_{n} T_{n} \biggr\| \\
&\ge& 1 + 2(\alpha-\eps) - \|T\|
~\ge~ 1+\|T\|-2\eps.
\eea

\begin{cor}\label{C6b}
If $X$ is a separable Banach space with the \DP, then $X$ does not
embed into an unconditional sum of reflexive spaces. In particular,
$X$ does not embed into a space with an unconditional basis.
\end{cor}

\Proof
Assume that $X$ embeds into an unconditional 
sum of Banach spaces $Y=\bigoplus_{n} X_{n}$ with associated
projections $P_{n}$ from $Y$ onto $X_{n}$. After replacing $X_{n}$ by
$P_{n}(X)$ we assume that the $X_{n}$ and $Y$ are separable. We may
also assume, by renorming $Y$ according to Theorem~\ref{T6a}, 
that the pair $(X,Y)$ has the \DP\ so
that by Theorem~\ref{T3} every weakly compact operator from $X$ into
$Y$ satisfies the Daugavet equation~(\ref{eq1.2}). Since the
embedding operator $J\dopu X\to Y$ has an expansion into a pointwise
unconditionally convergent series $J=\sum_{n=1}^{\infty}
\rest{P_{n}}{X}$, we deduce from Lemma~\ref{L6aa} that
$\rest{P_{n_{0}}}{X}$ is not weakly compact for some $n_{0}$, and
$X_{n_{0}}$ is not reflexive.
\eop
\bigskip

Now we use Lemma~\ref{L1} to produce $\ell_{1}$-copies in spaces with
the \DP. First, an extension of that lemma.

\begin{lemma}\label{L4}
If $(X,Y)$ has the \DP, then for every finite-dimensional subspace
$Y_{0}$ of $Y$, every $\eps_{0}>0$ and every slice $S(x^*_{0},\eps_{0})$
of $B_{X}$ there is a slice $S(x^*_{1},\eps_{1})$
of $B_{X}$ such that
\beq\label{eq3}
\|y+tx\|\ge (1-\eps_{0})(\|y\|+|t|) \qquad\forall y\in Y_{0},\ x\in
S(x^*_{1},\eps_{1}).
\eeq
\end{lemma}

\Proof
Let $\delta= \eps_{0}/2$ and pick a finite $\delta$-net
$\{y_{1},\ldots,y_{n}\}$ in $S_{Y_{0}}$. By a repeated application of
Lemma~\ref{L1}(a) we obtain a sequence of slices $S(x^*_{0},\eps_{0})
\supset S(x^{*(1)}, \eps^{(1)}) \supset \ldots \supset 
S(x^{*(n)}, \eps^{(n)})$ such that one has
\beq\label{eq4}
\|y_{k}+x \| \ge 2-\delta 
\eeq
for all $x\in S(x^{*(k)}, \eps^{(k)})$. Put $x^*_{1}= x^{*(n)}$ and
$\eps_{1}= \eps ^{(n)}$; then (\ref{eq4}) is valid for every $x\in
S(x^*_{1},\eps_{1})$ and $k=1,\ldots,n$. This implies that for every
$x\in S(x^*_{1},\eps_{1})$ and every $y\in S_{Y_{0}}$ the condition
$$
\|y+x\| \ge 2-2\delta = 2-\eps_{0}
$$
holds.

Let $0\le t_{1},t_{2}\le 1$ with $t_{1}+t_{2}=1$. If $t_{1}\ge
t_{2}$, we have for $x$ and $y$ as above
\bea
\|t_{1}x + t_{2}y\| &=&
\|t_{1}(x+y) + (t_{2}-t_{1})y\| \\
&\ge&
t_{1}\|x+y\| - |t_{2}-t_{1}|\, \|y\| \\
&\ge&
t_{1}(2-\eps_{0}) + t_{2}-t_{1} \\
&=&
t_{1}+t_{2}-t_{1}\eps_{0} ~\ge~ 1-\eps_{0},
\eea
and an analogous argument shows this estimate in case $t_{1}<t_{2}$.

This implies (\ref{eq3}), by the homogeneity of the norm and the
symmetry of $S_{Y_{0}}$.
\eop

\begin{theo}\label{T5}
If $X$ has the \DP, then $X$ contains a copy of $\ell_{1}$.
\end{theo}

\Proof
Using Lemma~\ref{L4} inductively, it is easy to construct a sequence
of vectors $e_{1},e_{2},\ldots$ and a sequence of slices
$S(x^*_{n},\eps_{n})$, $\eps_{n}=4^{-n}$, $n\in\N$, such that 
$e_{n+1}\in S(x^*_{n+1},\eps_{n+1}) $ and every
element of $S(x^*_{n+1},\eps_{n+1})$ is -- ``up to $\eps_{n}$'' --
$\ell_{1}$-orthogonal to $\lin\{e_{1},\ldots,e_{n}\}$, which means
$$
\|y+x\|\ge (1-\eps_{n})(\|y\|+\|x\|) \qquad\forall y\in
\lin\{e_{1},\ldots,e_{n}\},\ x\in S(x_{n+1}^{*},\eps_{n+1}).
$$
The sequence $(e_{n})$ is then equivalent to the unit vector basis in
$\ell_{1}$.
\eop

\bigskip
The proof even shows the stronger result that $X$ contains
asymptotically isometric copies  of $\ell_{1}$ in the sense of
\cite{DowJLT}.

Next, we study to what extent the \DP\  is hereditary.

Recall that an $M$-ideal in a Banach space $X$ is a closed subspace
$J$ such that $X^{*}$ decomposes as $X^{*}= V \oplus_{1} J^{\bot}$
for some closed subspace $V$ of $X^{*}$, where $J^{\bot}=\{ x^*\in
X^{*}\dopu \rest{x^*}{J}=0\}$.  Then $\{\rest{x^*}{J}\dopu x^*\in V\}$ is
linearly isometric to $J^{*}$, and we shall write
\beq\label{eq4a}
X^{*} = J^{*} \oplus_{1} J^{\bot}.
\eeq

\begin{prop}\label{P7}
The \DP\ is inherited by $M$-ideals.
\end{prop}

\Proof
Suppose $J$ is an $M$-ideal in a Banach space $X$ with the \DP. 
Let $y\in S_{J}$
and $\eps>0$, and let $x^*\in J^{*}\subset X^{*}$ with $\|x^*\|=1$.
Consider the slices
\bea
S_{1}&=& 
\{\xi\in J\dopu \|\xi\|\le 1,\ x^*(\xi)\ge 1-\eps\}, \\
S&=&
\{\xi\in X\dopu \|\xi\|\le 1,\ x^*(\xi)\ge 1-\eps/3\}.
\eea
By Lemma~\ref{L2}, there is some $x\in S$  such that $\|x+y\| \ge
2-\eps/3$; hence there is some $y^*\in S_{X^{*}}$ with $y^*(x+y)\ge
2-\eps/3$. Decompose $y^*=y^*_{1}+y^*_{2} \in J^{*}\oplus_{1} J^{\bot}$ so
that $1=\|y^*\| = \|y^*_{1}\| + \|y^*_{2}\|$. Therefore we have
$$
y^*(x) + y^*_{1}(y) \ge 2-\eps/3 
$$
so that $y^*(x)\ge 1-\eps/3$ and $y^*_{1}(y)\ge 1-\eps/3$.
Consequently, $\|y^*_{1}\|\ge 1-\eps/3$ and thus $\|y^*_{2}\|\le \eps/3$.

Let us now consider the $\sigma(X,J^{*})$-topology on $X$. An
application of the Hahn-Banach theorem shows that $B_{J}$ is
$\sigma(X,J^{*})$-dense in $B_{X}$ \cite[Remark~I.1.13]{HWW}. 
We may therefore find some
$\xi\in B_{J}$ satisfying $|y^*_{1}(\xi-x)| \le \eps/3$  and
$|x^*(\xi-x)|\le \eps/3$, i.e., $\xi\in S_{1}$, and we have
\bea
\|\xi+y\|
&=&
y^*(\xi+y) ~=~ y^*_{1}(\xi) + y^*_{1}(y) \\
&\ge&
y^*_{1}(x) + y^*_{1}(y) -\eps/3 \\
&\ge&
y^*_{1}(x) + y^*_{2}(x) + y^*_{1}(y) - 2\eps/3 \\
&& \qquad\mbox{since $\|y^*_{2}\|\le\eps/3$}  \\
&=&
y^*(x) + y^*_{1}(y) -2\eps/3  \\
&\ge&
2-\eps.
\eea
An application of Lemma~\ref{L2} completes the proof of the
proposition.
\eop
\bigskip

Obviously, if $X$ has the \DP\ and $J\subset X$ is an $M$-ideal,
then $X/J$ need not have the \DP; for example, if $X=C[0,1]$ and 
$J=\{f\in X\dopu f(0)=0\}$, then $X/J$  is one-dimensional and thus
fails the \DP.

We now prove a converse to Proposition~\ref{P7}, which can be
regarded as a version of the three space property for the \DP\ under
strong geometric assumptions.

\begin{prop}\label{7a}
If $J$ is an $M$-ideal in $X$ such that $J$ and $X/J$ share the \DP,
then so does $X$.
\end{prop}

\Proof
Suppose that $y\in S_{X}$, $x^{*}\in S_{X^{*}}$ and $\eps>0$ are
given as in Lem\-ma~\ref{L2}(ii). We decompose
$$
x^{*}= x_{1}^{*} + x_{2}^{*} \in J^{*}\oplus J^{\bot},\quad
\|x^{*}\| = \|x_{1}^{*}\| + \|x_{2}^{*}\|
$$
and from (\ref{eq4a}) we deduce that
$$
\|y\|=\max\biggl\{
\sup_{y^{*}\in B_{J^{*}}} |y^{*}(y)| ,
\sup_{y^{*}\in B_{J^{\bot}}} |y^{*}(y)| 
  \biggr\}  =1 .
$$

We shall first assume that
$$
\|[y]\|_{X/J}  = 
\sup_{y^{*}\in B_{J^{*}}} |y^{*}(y)| 
=1.
$$
Since $X/J$ has the \DP\ and $(X/J)^{*}=J^{\bot}$, there is some
$x_{0}\in X$ satisfying
$$
\|[x_{0}]\|=1, \quad
x_{2}^{*}(x_{0}) \ge (1-\eps)\|x_{2}^{*}\|, \quad
\|[x_{0}+y]\|\ge 2-\eps.
$$
Next, pick $\xi\in B_{J}$ with
\beq\label{eqM2}
x_{1}^{*}(\xi) \ge (1-\eps) \|x_{1}^{*}\|
\eeq
and use the 2-ball property of $M$-ideals \cite[Th.~I.2.2]{HWW} to
find some $\eta\in J$ with 
\beq\label{eqM3}
\|x_{0}\pm\xi -\eta\| \le 1+\eps.
\eeq
Obviously, $x:= x_{0}+\xi-\eta$ has the properties 
\bea
\|x\| &\le& 1+\eps, \\
x_{2}^{*}(x) &=& x_{2}^{*}(x_{0}) ~\ge~ (1-\eps)\|x_{2}^{*}\|, \\
\|x+y\| &\ge& \|[x+y]\|  ~=~  \|[x_{0}+y]\| ~\ge~ 2-\eps,
\eea
and it is left to estimate $x_{1}^{*}(x)$. Now we get from
(\ref{eqM3})
$$
 |x_{1}^{*}(\xi) \pm x_{1}^{*}(x_{0}-\eta)| \le
(1+\eps)\|x_{1}^{*}\|
$$
and hence from (\ref{eqM2}) 
$$
|x_{1}^{*}(x_{0}-\eta)| \le 2\eps \|x_{1}^{*}\|
$$
so that
$$
x_{1}^{*}(x) \ge (1-3\eps)\|x_{1}^{*}\|
$$
and finally
$$
x^{*}(x) \ge (1-3\eps)\|x_{1}^{*}\| + (1-\eps)\|x_{2}^{*}\|  
\ge 1-3\eps.
$$
After scaling $x$ appropriately we obtain (ii) of Lemma~\ref{L2}.

In the second part of the proof we suppose that
$$
\theta :=
\sup_{y^{*}\in B_{J^{\bot}}} |y^{*}(y)| <
\sup_{y^{*}\in B_{J^{*}}} |y^{*}(y)| 
=1.
$$
We shall need the following claim: There is some $\xi\in S_{J}$ such
that $\xi^{*}(y) \ge 1-3\eps$ whenever $\xi^{*}\in S_{J^{*}}$ and
$\xi^{*}(\xi)\ge 1-\eps$. In fact, we have a decomposition
$X^{**}=J^{\bot\bot}\oplus_{\infty} J^{*\bot}$ of the bidual space;
denote the projection from $X^{**}$ onto $ J^{\bot\bot}$ by $Q$. Now,
$$
1=\|y\|= \max\{ \|Qy\|, \|y-Qy\| \} = \max\{ \|Qy\|,\theta \} 
$$
and thus $\|Qy \|=1$. By the principle of local reflexivity, in the
version of \cite{BehLR2}, there is a linear operator $L\dopu
\lin\{y,Qy\}\to X$ such that
$\xi:=L(Qy)\in S_{J}$, $Ly=y$ and $\|L\|\le1+\eps$;
the point here is that $L$ maps $Qy\in  J^{\bot\bot}$ into $J$.
Clearly $\xi=\frac12 y + \frac12 (2\xi-y)$  and 
$$
\|2\xi-y\| = \|L(2Qy - y)\| \le (1+\eps) \|2Qy-y\| = 1+\eps.
$$
Hence, if $\xi^{*}\in S_{J^{*}}$, then $\xi^{*}(y)\le1$ and
$\xi^{*}(2\xi-y)\le 1+\eps$. Consequently $\xi^{*}(y)\ge1-3\eps$
whenever $\xi^{*}(\xi)\ge1-\eps$.

By assumption on $J$ and Lemma~\ref{L2} there is
some $x_{0}\in J$ such that
\beq\label{eqM5}
\|x_{0}\|=1 , \quad
x_{1}^{*}(x_{0}) \ge (1-\eps)\|x_{1}^{*}\|, \quad
\|x_{0}+\xi\|\ge 2-\eps.
\eeq
Next, pick $z\in B_{X}$ and $\xi_{0}^{*}\in S_{J^{*}}$ with the properties
$$
x_{2}^{*}(z)\ge (1-\eps)\|x_{2}^{*}\|, \quad
\xi_{0}^{*}(x_{0}+\xi) \ge 2-\eps
$$
so that 
\beq\label{eqM6}
\xi_{0}^{*}(x_{0})\ge 1-\eps, \quad \xi_{0}^{*}(\xi)\ge 1-\eps.
\eeq
By construction of $\xi$ we therefore have $\xi_{0}^{*}(y)\ge
1-3\eps$.
Using the 2-ball property of $M$-ideals again, we may find some
$\eta\in J$ with 
$$
\|z\pm x_{0}-\eta\| \le 1+\eps,
$$
and we let $x:= z+x_{0}-\eta$. As in the first part of the proof we
obtain from (\ref{eqM5}) and (\ref{eqM6})
$$
|x_{1}^{*}(z-\eta)| \le 2\eps\|x_{1}^{*}\|, \quad
|\xi_{0}^{*}(z-\eta)|\le 2\eps
$$
and thus 
\bea
x_{1}^{*}(x) &\ge& 
x_{1}^{*}(x_{0}) -2\eps\|x_{1}^{*}\| ~\ge~ (1-3\eps)\|x_{1}^{*}\|, \\
x_{2}^{*}(x)   &=&
x_{2}^{*}(z)   ~\ge~ (1-\eps) \|x_{2}^{*}\|
\eea
so that
\bea
\|x\| &\le& 1+\eps, \\
x^{*}(x) &\ge& 1-3\eps, 
\eea
and
$$
\|x+y\| ~\ge~ \xi_{0}^{*}(x+y) ~=~ 
\xi_{0}^{*}(z-\eta) + \xi_{0}^{*}(x_{0}+y) ~\ge~ 2-6\eps.
$$
Again, we see that (ii) of Lemma~\ref{L2} is fulfilled.
\eop
\bigskip

We next aim at proving that subspaces of Daugavet spaces having a
separable annihilator (such subspaces can be regarded as ``large'')
have the \DP\ as well. We shall use local $\ell_{1}$-complements as explained
in the following lemma.

\begin{lemma}\label{L8a}
If $X$ has the \DP, then for every $x\in S_{X}$, every $\eps>0$ and
every separable subspace $V$ of $X^{*}$ there is an element $x^{*}\in
S_{X^{*}}$ such that $x^{*}(x)\ge 1-\eps$  and
$\|x^{*}+v^{*}\|=1+\|v^{*}\|$ for all $v^{*}\in V$.
\end{lemma}

\Proof
Take a dense sequence $(v_{n}^{*})$ in $ S_{V}$ and a sequence of
weak$^{*}$ compact slices $S(x,\eps) \supset S(x_{1},\eps_{1})
\supset S(x_{2},\eps_{2}) \supset \ldots$  such that 
$$
\|x^{*}+v_{k}^{*}\|\ge 2-\frac1n \qquad\forall x^{*}\in
S(x_{n},\eps_{n}),\ k=1,\ldots,n;
$$
this is possible by a repeated application of
Lemma~\ref{L1}(b), as in the proof of Lemma~\ref{L4}.
Clearly, any $x^{*}\in \bigcap _{n\in\N} S(x_{n},\eps_{n})$ works.
\eop

\begin{cor}\label{C8b}
If $X$ has the \DP, then $X^{*}$ contains an isometric copy of
$\ell_{1}$. In particular, $X^{*}$ is neither strictly convex nor
smooth.
\end{cor}

\Proof
An isometric $\ell_{1}$-copy can be produced from Lemma~\ref{L8a} by
an obvious inductive procedure.
\eop

\begin{theo}\label{T8c}
Suppose $X$ has the \DP\ and $Y\subset X$ is a subspace with a separable
annihilator $Y^{\bot}$. Then $Y$ has the \DP.
\end{theo}

\Proof
Fix $y\in S_{Y}$ and $\eps>0$ and consider the slice
$$
S=\{y^{*}\in Y^{*} = X^{*}/Y^{\bot} \dopu \|y^{*}\|\le 1,\ y^{*}(y)\ge
1-\eps\}.
$$
Also, fix an element $[x_{1}^{*}]\in S_{X^{*}/Y^{\bot}}$. To prove
the theorem it suffices, by Lemma~\ref{L2}(iii), to find some
$[x_{2}^{*}]\in S$ such that $\|[x_{1}^{*}+x_{2}^{*}]\|=2$. To
achieve this, apply Lemma~\ref{L8a} with $x=y$ and
$V=\lin(\{x_{1}^{*}\}\cup Y^{\bot})$. We get a functional
$x_{2}^{*}\in S_{X^{*}}$ such that $x_{2}^{*}(y)\ge1-\eps$ and
$$
\|x_{2}^{*}+v^{*}\| = 1+\|v^{*}\| \qquad\forall v^{*}\in V.
$$
Then $[x_{2}^{*}]\in S$ and
\bea
\|[x_{1}^{*}+x_{2}^{*}]\| &=&
\inf \{ \|x_{1}^{*}+x_{2}^{*}+z^{*}\|\dopu z^{*}\in Y^{\bot} \} \\
&=&
\inf \{ 1 + \|x_{1}^{*}+z^{*}\|\dopu z^{*}\in Y^{\bot} \} \\
&& \qquad\mbox{(since $x_{1}^{*}+z^{*}\in V$)} \\
&=&
1+\|[x_{1}^{*}]\| ~=~ 2.
\eea
This completes the proof of the theorem.
\eop

\bigskip
In \cite{Kadets} the following property was investigated: $X$ has the
hereditary \DP\ if every finite-codimensional subspace  of $X$ has
the ordinary \DP. The preceding theorem shows that the hereditary \DP\
coincides with the usual one.

At the end of this section we study sums of pairs with the \DP.

\begin{lemma}\label{L8}
If $(X_{1},Y_{1})$ and $(X_{2},Y_{2})$ have the \DP, 
then so do
$(X_{1}\oplus_{1}X_{2}, Y_{1}\oplus_{1}Y_{2})$  and
$(X_{1}\oplus_{\infty}X_{2}, Y_{1}\oplus_{\infty}Y_{2})$.
\end{lemma}

\Proof
We first deal with $(X_{1}\oplus_{\infty}X_{2}, Y_{1}\oplus_{\infty}Y_{2})$. 
Let us consider 
$x^*_{j}\in X_{j}^{*}$, $y_{j}\in Y_{j}$ ($j=1,2$) with
$\|(y_{1},y_{2})\|=\max\{\|y_{1}\|,\|y_{2}\|\} =1$,
$\|(x^*_{1},x^*_{2})\| = \|x^*_{1}\| + \|x^*_{2}\| =1$. 
Assume without loss of generality that $\|y_{1}\|=1$. By
Lemma~\ref{L2} there is, given $\eps>0$, some $x_{1}\in X_{1}$
satisfying
$$
\|x_{1}\|=1,\quad x^*_{1}(x_{1})\ge \|x^*_{1}\|(1-\eps) , \quad
\|x_{1}+y_{1}\|\ge 2-\eps.
$$
Also, pick $x_{2}\in X_{2}$ such that
$$
\|x_{2}\|=1,\quad x^*_{2}(x_{2})\ge \|x^*_{2}\|(1-\eps).
$$
Then $\|(x_{1},x_{2})\|=1$, $\langle (x^*_{1},x^*_{2}), (x_{1},x_{2})
\rangle \ge 1-\eps$ and
$$
\|(x_{1},x_{2}) + (y_{1},y_{2}) \| \ge
\|x_{1} + y_{1}\| \ge 2-\eps.
$$
Thus, $(X_{1}\oplus_{\infty}X_{2}, Y_{1}\oplus_{\infty}Y_{2})$ has the \DP.

A similar calculation, based on Lemma~\ref{L2}(iii), shows that
$(X_{1}\oplus _{1}X_{2},\allowbreak Y_{1}\oplus_{1}Y_{2})$  has the \DP.
\eop

\bigskip
We remark that the converse of the above lemma is valid, too.

\begin{prop}\label{P9}
Suppose that $(X_{1},Y_{1}), (X_{2},Y_{2}), \ldots$ are pairs of
Banach spaces with the \DP.
Then $(c_{0}(X_{j}), c_{0}(Y_{j}))$ and $(\ell_{1}(X_{j}),
\ell_{1}(Y_{j}))$ have the \DP.
\end{prop}

\Proof
It follows from Lemma~\ref{L8} that $(X_{1}\oplus_{\infty}\cdots
\oplus_{\infty}X_{n}, Y_{1}\oplus_{\infty}\cdots
\oplus_{\infty}Y_{n})$, resp.\ $(X_{1}\oplus_{1}\cdots
\oplus_{1}X_{n}, Y_{1}\oplus_{1}\cdots
\oplus_{1}Y_{n})$, have the \DP\ for each $n\in \N$. Since the union
of these spaces is dense in $c_{0}(X_{j})$, $c_{0}(Y_{j})$,
$\ell_{1}(X_{j})$ or $\ell_{1}(Y_{j})$,
respectively, the result follows.
\eop

\bigskip
For $X_{j}=Y_{j}$
these results were first proved by Wojtaszczyk \cite{Woj92}
and, in a special case, by Abramovich \cite{Abra2} using different
approachs.

\mysec{3}{Quotients of $L_{1}$}%
In this section,
we consider the space $L_{1}=L_{1}[0,1]$. The Lebesgue measure is
denoted by $\mu$, and $\Sigma$ stands for the Borel $\sigma$-algebra
on $[0,1]$.

The following two properties which a subspace $X$ of $L_{1}$ might or
might not have will turn out to be relevant for the \DP\ of
$L_{1}/X$.
\medskip
\begsta
\item[\rm (I)]
For  every Borel set $A\subset [0,1]$ and
every $\eps>0$  there is a positive function $f\in L_{1}(A)$ such
that $\|f\|=1$ and
\beq\label{eq4.0}
(1-\eps) (|\lambda| + \|h\|)  
 \le \|\lambda f + h\| 
\le |\lambda| + \|h\|  
\eeq
for all $h\in X$, $\lambda\in\R$.
\item[\rm (II)]
For every $f_{0}\in L_{1}$, every Borel set $A\subset [0,1]$ and
every $\eps>0$  there is a positive function $f\in L_{1}(A)$ such
that $\|f\|=1$ and (\ref{eq4.0}) holds
for all $h\in \lin(X\cup\{f_{0}\})$, $\lambda\in\R$.
\endsta

\begin{lemma}\label{L4.1}
Let $X$ be a subspace of $L_{1}$.
\begsta
\item
{\rm (I)} implies\/  {\rm (II)}.
\item
If $X$ satisfies\/ {\rm (II)}, then $L_{1}/X$  has the \DP.
\endsta
\end{lemma}

\Proof
(a)
Let $f_{0}\notin X$, and let $P$ denote the projection from
$\lin(X\cup\{f_{0}\})$ onto $\lin\{f_{0}\}$ along $X$. Given
$\eps>0$, there is some $\delta>0$ such that
$$
\mu(B)<\delta \quad\Rightarrow \quad \|\chi_{B}f_{0}\|\le
\frac{\eps}{4\|P\|}.
$$
For $A$ as in (II), let $B\subset A$ be a Borel set of
measure ${<\delta}$. An application of (I) with $B$ in place of
$A$ and $\eps_{1}= \eps/(4+2\|P\|)$ in place of $\eps$ provides
us with a positive function $f\in L_{1}(B)$, $\|f\|=1$, such that
$$
\|\alpha f+g\| \ge (1-\eps_{1})(|\alpha|+\|g\|) \qquad\forall
\alpha\in\R,\ g\in X .
$$
Now let $h=\lambda f_{0}+g $ be a function in $\lin(X\cup \{f_{0}\})$
with $\|h\|=1$; i.e., $Ph= \lambda f_{0}$ and $|\lambda|\le \|P\|$.
To obtain (\ref{eq4.0}), it is enough to show that $\|f+h\|\ge
2-\eps$; see the last part of the proof of Lemma~\ref{L4}. In fact, we have
\bea
2-\|f+h\| &=&
\|f\| + \|h\| - \|f+h\| \\
&=&
\int_{B} \bigl( |f| + |\lambda f_{0}+g| - |f+\lambda f_{0}+g| \bigr)
                          \,d\mu \\
&& \qquad\mbox{(since $f=0$ on $[0,1]\setminus B$)} \\
&\le&
\int_{B} \bigl( |f| + |g| - |f+g| \bigr) +2|\lambda|
                   \int_{B}|f_{0}|\,d\mu \\
&\le&
 \bigl( \|f\| + \|g\| - \|f+g\| \bigr)
              +2|\lambda| \frac{\eps}{4\|P\|}  \\
&\le&
\eps_{1}(1+\|g\|) + \frac\eps2 \\
&\le&
(1+\|\Id-P\|)\|h\| \eps_{1} + \frac\eps2 \\
&\le&
\eps.
\eea

(b)
We shall verify condition~(ii) of Lemma~\ref{L2}. So let $[f_{0}]\in
S_{L_{1}/X}$, $F\in S_{X^{\bot}}\subset S_{L_\infty}$ and $\eps>0$ be
given; we are going to find an equivalence class $[f_{1}]$ in the slice
$S(F,\eps)$ such that $\|[f_{0}+f_{1}]\|\ge 2-2\eps$.

There is a Borel set $A\subset [0,1]$ such that
\beq\label{eq4.1}
1-\eps\le \rest{F}{A} \le 1 \ {\rm a.e.}
\quad {\rm or} \quad
-1\le \rest{F}{A} \le -1+\eps \ {\rm a.e.}
\eeq
Let us assume the former case. If $f$ satisfies the condition in
(II), then we have
$$
\langle F,f \rangle = \int_{A} F(t)f(t)\,dt \ge 1-\eps,
$$
i.e., $[f]\in S(F,\eps)$. If we put $f_{1}=f$, we get
\bea
\|[f_{0}+f_{1}]\| &=&
\inf \{ \|f_{1}+f_{0}+g\|\dopu g\in X\} \\
&\ge&
\inf \{ (1-\eps)(1+\|f_{0}+g\|)\dopu g\in X\} \\
&=&
(1-\eps) (1+\|[f_{0}]\|) \\
&=&
2-2\eps.
\eea

In the latter case of (\ref{eq4.1})  we put $f_{1}=-f$ to
obtain the same conclusion.
\eop

\begin{prop}\label{P4.2}
If $X$ is a reflexive subspace of $L_{1}$, then $L_{1}/X$ has the
\DP.
\end{prop}

\Proof
We shall use the previous lemma and stick to its notation. 
Now $B_{X}$ is uniformly integrable since it is
weakly compact \cite[p.~162]{Beau}. 
Thus, given $\eps>0$ there is some $\delta >0$ such
that
$$
\mu(B)\le \delta \qquad\Rightarrow\qquad \int_{B}|h(t)|\,dt \le \eps/2
\quad \forall h\in B_{X}.
$$
For a Borel set $A$ pick any Borel subset $B\subset A$ of
measure ${\le\delta}$, and let $f=\chi_{B}/\mu(B)$. Then we have for all
$h\in B_{X}$ and $|\lambda|+\|h\|=1$
\bea
\|\lambda f + h\| &=&
\int_{B} |\lambda f(t) + h(t)| \,dt  +  
\int_{[0,1]\setminus B}  |h(t)|\,dt  \\
&\ge&
|\lambda| -2 \int_{B} |h(t)|\,dt  + \int_{0}^{1} |h(t)|\,dt \\
&\ge&
1-\eps,
\eea
therefore property~(I) is fulfilled.
\eop
\bigskip

Our main objective in this section is to give an example of a Banach
space with the \DP\ not containing a copy of $L_{1}$, thus showing
that a conceivable generalization of Theorem~\ref{T5} is not valid.
Our example will be the quotient space $L_{1}/Y$, with $Y$ a space
constructed by Talagrand \cite{Tala5} as a counterexample to the
three space problem for $L_{1}$. Before we present several crucial
features of $Y$, we formulate a lemma.

\begin{lemma}\label{L4.3}
Let $g_{m}= f_{m}+h_{m}\in L_{1}$ be a sequence of bounded functions
on $[0,1]$ such that $\sup_{m} \|f_{m}\|_{L_{\infty}} =:M_{1}
<\infty$. Assume also that the supports $\Delta_{m}$ of $h_{m}$  are
small in the sense that $\sum_{m=1}^{\infty} \mu(\Delta_{m}) <\infty$
and that the sequence $(g_{m})$ is equivalent to the standard
$\ell_{1}$-basis, i.e.,
$$
\biggl\| \sum_{m=1}^{\infty} \lambda_{m} g_{m} \biggr\|
\le \sum_{m=1}^{\infty} |\lambda_{m}| \le
M_{2} \biggl\| \sum_{m=1}^{\infty} \lambda_{m} g_{m} \biggr\|
$$
for some $M_{2}<\infty$ and all $(\lambda_{m})\in \ell_{1}$.
Then for every $\eps>0$ there is some $\delta\in (0,\eps/2)$ such
that whenever $A\in \Sigma$, $\mu(A)>\eps$, there is a subset
$B\subset A$, $\mu(B)=\delta$, satisfying
\beq\label{eqT1}
\|\chi_{B}f\| \le \eps \|f\| \qquad\forall f\in X:= \linq\{g_{m}\dopu
m\in \N\}.
\eeq
\end{lemma}

\Proof
Select a number $N\in\N$ for which $\sum_{m>N}\mu(\Delta_{m}) <
\eps/2$. Denote $M_{3}:=\max_{k\le N} \|h_{k}\|_{L_{\infty}}$ and put 
$\delta= \eps/ (2+M_{2}(M_{1}+M_{3}))$. Then for every $A\in \Sigma$,
$\mu(A)>\eps$, there is a subset $B\subset A$, $\mu(B)=\delta$, such
that $B\cap \bigl( \bigcup_{m>N} \Delta_{m} \bigr) =\emptyset$, since
$\delta\le\eps/2$. Take
any $f\in S_{X}$, $f=\sum_{m=1}^{\infty}\lambda _{m} g_{m}$. Then
$\sum_{m=1}^{\infty}|\lambda_{m}| \le M_{2}$, and we obtain 
$$
\|\chi_{B}f\| = \int_{B} \biggl| \sum_{m=1}^{\infty} \lambda_{m}f_{m}
   + \sum_{m=1}^{N} \lambda_{m}h_{m} \biggr| \,d\mu
\le \mu(B) M_{2} (M_{1}+M_{3}) <\eps,
$$
which proves (\ref{eqT1})
\eop
\bigskip

Let us now describe the structure of Talagrand's example. 
In \cite{Tala5} he
constructs a space generated by a double sequence
$g_{m,n}=f_{m,n}+h_{m,n}$ such that for every fixed $n$ the sequence
$(g_{m,n})_{m}$ meets the conditions of Lemma~\ref{L4.3}. Moreover,
if one denotes $X_{n}= \linq\{g_{m,n}\dopu m\in\N\}\subset L_{1}$,
then each $X_{n}$ is isomorphic to $\ell_{1}$, with Banach-Mazur
distance tending to $\infty$, though, and $\linq \bigcup_{n}X_{n}$ is
isomorphic to the $\ell_{1}$-sum of the $X_{n}$ with isomorphism
constant ${\le20}$, and not only each $X_{n}$ but every finite sum of
these spaces meets the conditions of Lemma~\ref{L4.3}. 
The $X_{n}$ are constructed so as to consist only
of very ``peaky'' functions: For every $\eps>0$ there is an
$n=n(\eps)$ such that \cite[Th.~3.1]{Tala5}
\beq\label{eqT2}
\mu(\{|f|\ge\eps\}) \le \eps \qquad\forall f\in S_{X_{n}}.
\eeq
Finally, whenever $(Y_{n})$ is a subsequence of $(X_{n})$, then
$L_{1}/\linq \bigcup_{n} Y_{n}$ fails to contain a copy of $L_{1}$
\cite[p.~26]{Tala5}.

\begin{theo}\label{T4.4}
There is a subspace $Y$ of $L_{1}$ such that $L_{1}/Y$ has the \DP,
but fails to contain a copy of $L_{1}$.
\end{theo}

\Proof
The space $Y$ will be the closed linear span of a certain subsequence
$(Y_{n})$ of Talagrand's sequence $(X_{n})$.  As reported above, it
is a deep result due to Talagrand that $L_{1}/Y$ does not contain a
copy of $L_{1}$.

In order to obtain the \DP\ for $L_{1}/Y$ by means of
Lem\-ma~\ref{L4.1}, we determine the $Y_{n}$ as follows. First, we let
$Y_{1}=X_{1}$.     Using inductively Lemma~\ref{L4.3} and condition
(\ref{eqT2}) we select a subsequence $(Y_{n})$ to fulfill the
following conditions:
\medskip
\begsta
\item
For every $n\in\N$ there is a $\delta_{n}\in (0,2^{-n-1})$ such that for
every $A\in\Sigma$, $\mu(A)>2^{-n}$, there is some $B\subset A$,
$\mu(B)=\delta_{n}$, such that 
\beq\label{eqT1a}
\|\chi_{B}f\| \le 2^{-n} \|f\| \qquad\forall f\in
Y_{1}\oplus\cdots\oplus Y_{n};
\eeq
in addition we pick $\delta_{n}<\delta_{n-1}$.
\item
For every $g\in S_{Y_{n+1}}$ we have
\beq\label{eqT2a}
\mu(\{ |g|\ge 10^{-n}\delta_{n} \}) \le 10^{-n}\delta_{n}.
\eeq
\endsta
Put now $Y=\linq \bigcup_{n}Y_{n}$ and let us prove the validity of
condition~(I) of Lemma~\ref{L4.1}. 

Denote by $M$ the smallest constant such that 
\beq\label{eqT3}
M\biggl\| \sum_{n=1}^{\infty} \lambda_{n}f_{n} \biggr\| \ge
\sum_{n=1}^{\infty} |\lambda_{n}| 
\eeq
for every choice of $f_{n}\in S_{Y_{n}}$ ($n\in\N$) and for all
$(\lambda_{n})\in \ell_{1}$; in fact, $M\le20$
\cite[Prop.~4.2]{Tala5}. Fix $\eps>0$ and a set $A\in\Sigma$,
$\mu(A)>\eps$. Pick $N\in\N$ such that
\beq\label{eqT4}
2^{-N+1} (M+2) \le \eps,
\eeq
and apply (a) to obtain some $B\in\Sig$ satisfying (\ref{eqT1a}) for
$n=N$. Put $g_{0}=\chi_{B}/\mu(B)$ and consider a function
$h=\sum_{n=1}^{\infty} \lambda_{n}f_{n}$ with $\|h\|=1$ where
$f_{n}\in S_{Y_{n}}$ for all $n$  and $(\lambda_{n})\in\ell_{1}$. We
need to prove that
\beq\label{eqT6}
\|g_{0}+h\|\ge 2-\eps.
\eeq

Denote
$$
h_{0}=\sum_{n=1}^{N} \lambda_{n}f_{n} , \quad
h_{1}= h-h_{0}, \quad
D=\{t\in B\dopu |h_{1}(t)|\le 2^{-N}\}.
$$
By (\ref{eqT3})  we know that 
$$
\sum_{n=1}^{\infty} |\lambda_{n}| \le M,
$$
so for every $t\in D$ at least one of the inequalities 
$$
|f_{n}(t)|\ge 10^{-n}, \quad n=N+1, N+2, \ldots
$$
has to be true, hence by (\ref{eqT2a}) 
\beq\label{eqT8}
\mu(D) \le \sum_{n>N} \mu(\{ |f_{n}|\ge 10^{-n} \}) \le
10^{-N}\delta_{N}.
\eeq
To prove (\ref{eqT6})  let us note that
\bea
2-\|g_{0}+h\|
&=&
\|g_{0}\| + \|h\| - \|g_{0}+h\| \\
&=&
\int_{B} \bigl( |g_{0}| + |h| - |g_{0}+h| \bigr) \,d\mu \\
&=&
\int_{B\setminus D} \bigl( |g_{0}| + |h| - |g_{0}+h| \bigr) \,d\mu \\
&& \qquad
  + \int_{D} \bigl( |g_{0}| + |h| - |g_{0}+h| \bigr) \,d\mu \\
&\le&
\int_{B\setminus D} 2|h|\,d\mu  +  \int_{D} 2|g_{0}|\,d\mu.
\eea
The last two integrals can be estimated from above as follows:
\bea
\int_{B\setminus D} |h|\, d\mu
&\le&
\int_{B} |h_{0}|\,d\mu + \int_{B\setminus D} |h_{1}|\,d\mu \\
&\le&
2^{-N} \sum_{n=1}^{N} |\lambda_{n}| + 2^{-N} \delta_{N} \\
&&  \qquad\mbox{(by (\ref{eqT1a}) and by definition of $D$)} \\
&\le & 
2^{-N} (M+1).
\eea
Further, by definition of $g_{0}$ and (\ref{eqT8})
$$
\int_{D} |g_{0}|\,d\mu  = \frac{\mu(D)}{\mu(B)} \le 10^{-N}
$$
so that 
$$
2-\|g_{0}+h\| \le 2^{-N} (M+2)\cdot2 \le \eps
$$
by (\ref{eqT4}), and the proof of (\ref{eqT6}) and thus the proof of
Theorem~\ref{T4.4} is completed.
\eop
\bigskip

The common feature of Proposition~\ref{P4.2} and Theorem~\ref{T4.4}
is that we factor by a subspace with the Radon-Nikod\'ym property.
Recall from \cite{Woj92} or \cite{Dirk10} that $L_{1}/H_{1}$ has the
\DP, too (since its dual does). So a natural question that has
remained open is whether a quotient $L_{1}/X$ by a subspace with the
Radon-Nikod\'ym property always has the \DP\ or can at least  be so
renormed.

\mysec{4}{The anti-Daugavet property}%
In this part we consider operators $T$ from a Banach space $X$ into
itself. It is an easy exercise to prove that 
\beq\label{eq3.0}
\|\Id+T\|=1+\|T\|
\eeq 
if $T$ is a
bounded linear operator for which $\|T\|$ belongs to the spectrum of
$T$. Therefore, we say that a Banach space $X$ has the {\em anti-\DP\/} for a
class $\cal M$ of operators if, for $T\in\cal M$, the equivalence
$$
\|\Id+T\|=1+\|T\|  \quad \Longleftrightarrow  \quad  \|T\|\in\sigma(T)
$$
holds. Again, it is enough to consider operators of norm~$1$. If
${\cal M}=L(X)$, we simply speak of the anti-\DP. 

We shall use some results on finite representability, ultrapowers and
superreflexivity which may be found in the monograph \cite{Beau}.
Geometric notions such as uniform convexity and smoothness are
discussed there, too; see also \cite{Die-LNM} or \cite{HHZ}.

It was proved in
\cite{AbraAB} that uniformly convex and uniformly smooth spaces share
the anti-\DP, and locally uniformly convex spaces have the anti-\DP\
for compact operators.
In \cite{Kadets} an equivalent geometric condition
(see Definition~\ref{D3.1}(a))
for the anti-\DP\ in finite-dimensional spaces was presented. 
We shall now extend these results 
and introduce some geometric properties of the unit sphere.

\begin{defi}\label{D3.1}  \mbox{ }
\begsta
\item
We say that a Banach space $X$ is {\em alternatively convex or smooth\/} 
(acs) if for all $x,y\in S_{X}$ and $x^{*}\in S_{X^{*}}$
the implication
\beq\label{eq3.1}
x^{*}(x)=1,\ \|x+y\|=2 \quad\Rightarrow\quad x^{*}(y)=1
\eeq
holds.
\item
We say that a Banach space $X$ is {\em locally uniformly
alternatively convex or smooth\/} 
(luacs) if for all $x_{n},y\in S_{X}$ and $x^{*}\in S_{X^{*}}$
the implication
\beq\label{eq3.2}
x^{*}(x_{n})\to1,\ \|x_{n}+y\|\to2 \quad\Rightarrow\quad x^{*}(y)=1
\eeq
holds.
\item
We say that a Banach space $X$ is {\em uniformly
alternatively convex or smooth\/} 
(uacs) if for all $x_{n},y_{n}\in S_{X}$ and $x_{n}^{*}\in S_{X^{*}}$
the implication
\beq\label{eq3.3}
x_{n}^{*}(x_{n})\to1,\ \|x_{n}+y_{n}\|\to2 \quad\Rightarrow\quad 
x_{n}^{*}(y_{n})\to1
\eeq
holds.
\endsta
\end{defi}

Geometrically, the (acs)-property means some smoothness of the norm
at points lying on a line segment in the unit sphere. Precisely, $X$
is (acs) if and only if whenever $\co\{x,y\}\subset S_{X}$, then $x$
and $y$ are smooth points of the unit ball of $\lin\{x,y\}$.

We remark that uniformly convex spaces and uniformly smooth spaces
are (uacs), and locally uniformly convex spaces are (luacs).
In \cite{Kadets} it was proved that in finite-dimensional spaces
(where by a compactness argument (acs), (luacs) and (uacs) are equivalent)
(acs) is necessary and sufficient for the anti-\DP.
We intend to characterize the anti-\DP\ for the class of compact
operators on infinite-dimensional spaces 
by means of the (luacs)-property. To this end, we shall
need a lemma.

\begin{lemma}\label{L3.2}
Suppose $X$ is\/ {\rm (acs)} and $T\dopu X\to X$ is a weakly compact
operator with $\|T\|=1$ and $\|\Id+T\|=2$.    Suppose in addition that
$\|x+Tx\|=2$ for some $x\in S_{X}$. Then $1$ is an eigenvalue of $T$.
\end{lemma}

\Proof
Consider a functional $x^{*}\in S_{X^{*}}$ such that $x^{*}(x)=1$. 
By the (acs)-property of $X$ one has $x^{*}(Tx)=1$. Therefore
$x_{1}^*:= T^{*}(x^{*})$ attains the value~$1$  at $x$ and hence
belongs to $S_{X^{*}}$. Again, using  (\ref{eq3.1}) we obtain
$x_{1}^{*}(Tx)=1$. Applying the same argument inductively shows that
$x^{*}(T^{n}x)=1$ for all $n\in \N$. This implies that
$$
K:= \coq \{T^{n}x\dopu n\in\N\} \subset  
\{v\in X\dopu x^{*}(v)=1\}  \cap B_{X}  \subset S_{X},
$$
in particular $0\notin K$. Also, $K$ is a weakly compact convex set,
since $\{T^{n}x\dopu n\ge1\} = T(\{T^{n}x\dopu n\ge0\}) $,
which is relatively weakly compact, and $T$ maps
$K$ into $K$. By the Schauder (or Markov-Kakutani) fixed point
theorem, $T$ has a fixed point in $K$, which is a non-zero
eigenvector for the eigenvalue~$1$.
\eop

\begin{theo}\label{T3.3}
For a Banach space $X$, the following conditions are equivalent:
\begaeq
\item
$X$ has the anti-\DP\ for compact operators.
\item
$X$ has the anti-\DP\ for operators of rank~$1$.
\item
$X$ is\/ {\rm (luacs)}.
\endaeq
\end{theo}

\Proof
(i) $\Rightarrow$ (ii) is evident. For the proof of (ii) $\Rightarrow$
(iii), assume that $X$ fails to be (luacs). Then there is a
functional $x^{*}\in S_{X^{*}}$ and there are elements $x_{n},y\in
S_{X}$ such that $x_{n},y\in S_{X}$  and $\|x+y_{n}\|\to 2$,
$x^{*}(x_{n})\to 1$, but $x^{*}(y)<1$. Consider the operator $T\dopu
X\to X$  defined by $Tv= x^{*}(v)y$. Then $\|T\|=1$ and
$\|\Id+T\|=2$, since
\bea
\|\Id + T\| &\ge&
\limsup \|x_{n}+Tx_{n}\| \\
&=&
\limsup \|x_{n}+ x^{*}(x_{n})y\| \\
&=&
\limsup \|x_{n}+y\| ~=~ 2 .
\eea
Thus $T$ satisfies (\ref{eq3.0}), but 
$1\notin \sigma(T)$ because of $x^{*}(y)<1$; so $X $ fails 
the anti-\DP\ for rank~$1$ operators.

(iii) $\Rightarrow$ (i):
Let $T$ be a compact operator with $\|T\|=1$ and $\|\Id+T\|=2$. Then
there is a sequence $(x_{n})\subset S_{X}$ for which
$\|Tx_{n}+x_{n}\|\to 2$. By compactness of $T$ we may assume that
$Tx_{n}\to y\in S_{X}$. Now consider $x^{*}\in S_{X^{*}}$ such that
$x^{*}(y)=1$; then $x^{*}(Tx_{n})\to 1$. Put $y^{*}= T^*x^{*}$; then
we have $\|y^{*}\|\le 1$, and from $y^{*}(x_{n})= x^{*}(Tx_{n})\to1$
we deduce that actually $\|y^{*}\|=1$. So we have $\|x_{n}+y\|\to2$
and $y^{*}(x_{n})\to 1$, and from (\ref{eq3.2}) we get that
$y^{*}(y)=1$. But now
$$
\|y+Ty\| \ge x^{*}(y+Ty) = x^{*}(y) + y^{*}(y) =2,
$$
so Lemma~\ref{L3.2} implies that $1\in \sigma(T)$; hence we
obtain~(i).
\eop
\bigskip

We now turn to the relation of the (uacs)-property and the anti-\DP.

\begin{lemma}\label{L3.4}
If $X$ is\/  {\rm (uacs)}, then $X$ is superreflexive.
\end{lemma}

\Proof
The (uacs)-property provides a uniform restriction on the structure
of 2-dimensional subspaces of $X$: For all $\eps>0$ there exists some
$\delta>0$ such that if $\|x+y\|>2-\delta$, $x,y\in S_{X}$, and
$x^{*}\in S_{(\lin\{x,y\})^{*}}$ with $x^{*}(x)>1-\delta$, then
$x^{*}(y)>1-\eps$. Therefore, not every 2-dimensional Banach space is
finitely representable in $X$, and by \cite{Kadets82} $X$ has to be
superreflexive.
(In fact, $X$ is uniformly non-square, which is enough to imply
superreflexivity by a theorem due to James; see \cite[p.~261]{Beau}.)
\eop

\begin{theo}\label{T3.5}
If $X$ is\/ {\rm (uacs)}, then $X$ has the anti-\DP.
\end{theo}

\Proof
Let $T\dopu X\to X$ be an operator of norm~$1$ such that
$\|\Id+T\|=2$. Consider an ultrapower $X^{\frak{U}}$ of $X$; $X^{\frak
U}$ is the factor space $\ell_{\infty}(X)/c_0^{\frak{U}}(X)$ with the
norm $\|[(x_{n})]\|=\lim_{\frak U}\|x_{n}\|$, where $c_0^{\frak U}(X)$
consists of those sequences in $X$ that tend to zero along the
ultrafilter~$\frak U$. Define $T^{\frak U}\dopu X^{\frak U}\to
X^{\frak U}$  by $T^{\frak U}[(x_{n})] = [(Tx_{n})]$. 

The main
advantage of considering $T^{\frak U}$ is that $\Id^{\frak
U}+T^{\frak U}$ attains its norm.
Indeed, if $x_{n}^{0}\in S_{X}$ are chosen so that $\|Tx_{n}^{0}\|\to 1$,
then $x^{0}:= [(x_{n}^{0})]\in S_{X^{\frak U}}$ and $\|T^{\frak
U}(x^{0})\|=1$. By Lemma~\ref{L3.4}, $X^{\frak U}$ is reflexive
and thus $T^{\frak U}$ is weakly compact. It is evident from the
definition that if $X$ is (uacs) and $Y$ is finitely representable in
$X$, then $Y$ is (uacs), i.e., (uacs) is a superproperty. Thus we
conclude that $X^{\frak U}$ is (acs), and we may apply
Lemma~\ref{L3.2}. So there is an eigenvector $x^{1}=[(x^{1}_{n})]\in
S_{X^{\frak U}}$ with $T^{\frak U}x^{1}=x^{1}$. This means that
$\|Tx^{1}_{n}-x^{1}_{n}\|$ tends to zero along $\frak U$, and 
$1$ is an approximate eigenvalue of $T$, in particular $1\in
\sigma(T)$.
\eop

\bigskip
As we said, it is evident that (uacs) is a superproperty. Also, the
super anti-\DP\ coincides with (uacs), but we doubt whether the
anti-\DP\ is a superproperty itself.

%
%
%
%
\typeout{References}

%
%
%
\small
\bigskip
\noindent
Faculty of Mechanics and Mathematics, Kharkov State University,
pl.~Svobody,~4, 310077~Kharkov, Ukraine

\smallskip\noindent
I.~Mathematisches Institut, Freie Universit\"at Berlin,
Arnimallee 2--6, \\
D-14\,195 Berlin, Germany; \ 
e-mail: 
werner@math.fu-berlin.de

\end{document}